\documentclass[11pt, leqno, amsfonts]{article}
\usepackage{amsmath,latexsym}

\renewcommand{\ref}{\par\noindent\hangindent 12pt}

\usepackage{amsmath}

\usepackage[dvips]{graphicx}

\makeatletter
\newcommand{\@makemycaption}[2]{%
\vspace{10pt}%
{\textbf{#1}:#2\par}%
}
\renewcommand{\figure}{%
\let\@makecaption\@makemycaption\@float{figure}}
\renewcommand{\table}{%
\let\@makecaption\@makemycaption\@float{table}}
\makeatother

\date{}

  {\end{list}%
  }%
\newcommand{\bx}{{\bf x}}
\newcommand{\by}{{\bf y}}
\newcommand{\bt}{{\bf t}}
\newcommand{\bu}{{\bf u}}
\newcommand{\bv}{{\bf v}}
\newcommand{\bw}{{\bf w}}
\newcommand{\bb}{{\bf b}}
\newcommand{\bn}{{\bf n}}

\newcommand{\bk}{{\bf k}}
\newcommand{\be}{{\bf e}}
\newcommand{\bs}{{\bf s}}
\newcommand{\noi}{\noindent}

\newcommand{\bX}{{\bf X}}
\newcommand{\bY}{{\bf Y}}
\newcommand{\bZ}{{\bf Z}}

\newcommand{\mZ}{{\bf Z}}
\newcommand{\mR}{{\bf R}}

\newcommand{\wtd}{\widetilde}
\newcommand{\cov}{{\rm Cov}}
\newcommand{\wQ}{\widetilde Q}

\begin{document}
\centerline{\bf MAXIMA OF MOVING SUMS IN A POISSON RANDOM FIELD} 

\vskip 0.5in

\centerline{ BY HOCK PENG CHAN}

\bigskip 
\centerline{\it National University of Singapore}

\vskip 0.5in

\centerline{\bf Abstract} 

\smallskip The extremal tail probabilities
of moving sums in a marked Poisson random field is 
examined here. These sums are computed by adding up
the weighted occurrences of events lying within a
scanning set of fixed shape and size. Change of measure and analysis of
local random fields are used to provide tail probabilities.
The asymptotic constants are initially expressed in a form that seems hard to evaluate
and do not seem to provide any additional information on the properties 
of the constants. A more sophisticated approach is then undertaken giving
rise to an expression that is not only neater but also able to provide computable
bounds. The technique used to obtain this constant can also be modified to work
on continuous processes. 

\vskip 2in

\noi Abbreviated Title: MAXIMA OF POISSON MOVING SUMS

\bigskip

\noi Supported by grants from the National 
University of Singapore.

\bigskip 
\noi AMS 2000 {\it subject classifications}. 
Primary; 60F10; secondary 60G10, 60G55.

\bigskip 
\noi {\it Key words and phrases}. Change of measure, large
deviations, marked Poisson process, moving sums, Poisson clumping, scan
statistics. 

\newpage 
{\bf 1. Introduction}. The maxima of moving
averages in Gaussian random fields in dimension
$d>1$ was studied in Siegmund and Worsley (1995) and 
Shafie, Sigal, Siegmund and Worsley (2003), with applications
in imaging and signal detection. Two key techniques
used are (i) the Karhunen-Lo\`{e}ve expansion with 
the volume of tube formula and (ii) the Euler 
characteristic; see Adler (2000) for an overview of
the research area and also Taylor, Takemura and Adler
(2005) and Taylor (2006) for more recent developments.

The maximum of moving sums in Poisson random fields,
more commonly known as scan statistics in the
statistical literature, also have widespread applications
in molecular biology, epidemiology, geostatistics
and image analysis, cf. Cressie (1993), 
Anderson and Titterington (1997), Glaz, Naus and
Wallenstein (2001) and Chan and Zhang (2007), but
the tail probability approximations are in
comparison not as
well developed for $d > 1$. While the tail probabilities
of these sums have been studied in Naus (1965), Loader (1991)
and Alm (1997), restrictions to rectangular scanning 
sets have been imposed for analytical convenience. 

We set out here to study the tail probabilities
of the maxima of moving sums with minimal restrictions
on the choice of scanning sets. A theory parallel to the study
of tail probabilities in Gaussian or Gaussian-like
random fields in the classical framework of Pickands (1969), Bickel and Rosenblatt
(1973), Qualls and Watanabe (1973), Piterbarg (1996)
and Chan and Lai (2006) is developed here.
We also consider a more general marked Poisson random
field, which is motivated by recent developments in
molecular biology, see for example Chan and Zhang (2007). 
This generalization entails careful consideration of overshoots
in special cases of scanning sets that is not 
required in Poisson random fields. Berman (1982) and Albin (1990)
have also studied tail probabilities of stationary
processes but their limiting results are of a
different type and do not apply here.

The first expression of the tail probability is 
stated in Theorem 1 in Section 3. Lemma 1 is the
basic building block of Theorem 1, providing the
extremal tail probability over a local domain by
using a change of measure approach. The expression
of this tail probability requires a description of
an induced local random field around the boundary
of the scanning set and this is provided in Section 2.
The technical details of how these building blocks can be
combined together to provide the tail probability
of the maxima of the sums over the whole domain,
via an adaption of the Pickands-Qualls-Watanabe
technique, is given in the Appendix. In Section 4,
we provide an alternative expression of 
the asymptotic constants in Theorem 1 via a more refined technique and obtain
bounds of these constants. In Section 5, we adapt this technique on continuous
valued random fields and show that it provide constants that looks like a differential
form of the constants obtained via the beautiful Poisson clumping heuristic shown
in Aldous (1989). Some bounds obtained from the new expression are surprisingly
accurate.  

\bigskip
{\bf 2. Definitions, notations and a local Poisson
random field}.
Let $D$ and $B$ be Jordan-measurable (bounded) subsets of ${\mR}^d$.
For vectors $\bt=(t_1,\ldots,t_d)$ and $\bu = (u_1,\ldots,u_d)$, we
shall use the notation $\bt \succ \bu$ to denote $t_j \geq u_j$ for all
$j$. We shall also let ${\bf 0} = (0,\ldots,0)$
and ${\bf 1}=(1,\ldots,1)$. Let $\sigma_k(\cdot)$,
$k = d$ or $k=d-1$ 
denote the $k$-dimensional volume of a $k$-dimensional manifold
in $\mR^d$. 
For any $A \subset \mR^d$, $b \in \mR$ and $c \in \mR^d$, we shall let
$\# A$ denote the number of elements in $A$ and $c+bA = \{ c+ba: a \in A \}$.
We shall also use $\| \cdot \|$ to denote $L_2$-norm
and $\| \cdot \|_\infty$ to denote $L_\infty$-norm.

Assume that the boundary $\partial B$
can be expressed as a finite union of 
smooth $(d-1)$-dimensional
submanifolds possibly with boundary (see Spivak (1965)
p113 for the definition). 
For example, if $B = \{ \bt: \| \bt \|_\infty \leq 1 \}$, a cube
of length 2, then $\partial B$ is a
union of $2d$ faces, each a smooth 
$(d-1)$-dimensional submanifold with boundary. 

Let $\bX = \{ (\bt_i, X_i): i \geq 1 \}$ be a marked
point process on ${\mR}^{d+1}$, 
characterized by $F$, a distribution
function of the marks $X_i$ and $\lambda > 0$, 
the rate of events occurring. Hence 
for any set $A$, Borel subset of ${\mR}^{d+1}$, 
$\# \{ i: (\bt_i, X_i) \in A \}$ follows a Poisson
distribution with mean $\lambda \int_A d \bt \times dF(x)$. 
Moreover, for any two
disjoint Borel sets $A$ and $C$, $\# \{ i: (\bt_i,X_i) 
\in A \}$ and $\# \{ i:
(\bt_i,X_i) \in C \}$ are independent random variables.
From the above description, 
we may assume without loss of generality that
$X_1,X_2,\ldots$ are independent and identically
distributed (i.i.d.) random variables having distribution $F$ and independent
of $\{ \bt_i: i \geq 1 \}$; a Poisson process with rate $\lambda$.

Let $\mu=EX_1$ and $M(\theta)=E e^{\theta X_1}$. 
Assume that $P \{ X_1 > 0 \} > 0$ and $\Theta := \{ \theta:
M(\theta) < \infty \}$ is an
open neighborhood of 0. For any set $A$, Borel subset
of ${\mR}^d$, define the sum $S(A) = 
\sum_{\bt_i \in A} X_i$. We analyze here the tail probability  
$$p_\lambda := P_\lambda \{ \sup_{\bv \in D} S(\bv+B) \geq \lambda c \}
\leqno(2.1)
$$
as $\lambda \rightarrow \infty$, 
for given $c>\max \{ 0, \mu \sigma_d(B) \}$. 

Through an appropriate transformation, we can also
look at the limiting probability of $p_\lambda$ as
one involving fixed Poisson rate $\lambda_0 > 0$
and increasingly large scanning sets. Let $g_\lambda
= (\lambda/\lambda_0)^{1/d}$ be the scaling constants. Then
$p_\lambda = P_{\lambda_0} \{ \sup_{\bv \in g_\lambda
D} S(\bv+g_\lambda B) \geq \lambda c \}$. For notational
simplicity, the analysis here looks at 
$p_\lambda$ in terms of (2.1) but the presence
of such transformations has important practical
implications. 

We will now proceed with the description of a limiting local
random field $\bY = \{ Y(\bu): \bu \in {\bf R}
\}$, that is derived from both the distribution 
$F$ and the geometry of the boundary $\partial
B$. For a given $c > \max \{ 0, \mu \sigma_d(B) \}$, let $\theta_c
>0$ and distribution $F_c$ satisfy
$$M'(\theta_c) = c/\sigma_d(B) \mbox{ and } F_c(dx)=
e^{\theta_c x} F(dx)/M(\theta_c), \leqno(2.2)
$$
where $'$ here denotes first derivative.

Let $\bZ^{(1)} = \{ (\bv_i^{(1)}, Z_i^{(1)}): i \geq 1 \}$
be a marked Poisson process such that $\{ \bv_i^{(1)}: i
\geq 1 \}$ is a Poisson process with rate 1 on the
domain $\partial B \times [0,\infty)$ and $Z_1^{(1)},
Z_2^{(1)},\ldots$ are i.i.d. with distribution $F$.
Let $\bZ^{(2)} = \{ (\bv_i^{(2)},Z_i^{(2)}):
i \geq 1 \}$ be a marked Poisson process independent
of $\bZ^{(1)}$, such that
$\{ \bv_i^{(2)}: i \geq 1 \}$ is a Poisson process
with rate $M(\theta_c)$ on the domain $\partial B
\times (-\infty,0)$ and $Z_1^{(2)},Z_2^{(2)},\ldots$
are i.i.d. with distribution $F_c$. 
Let $\bn_\bt$ be the unit normal vector of $\bt
\in \partial B$ away from $B$ and let $\cdot$ denote dot
product. For $\bu \in {\bf R}^d$, let
$$\begin{aligned} Y^{(j)}(\bu) = & \sum_{i: \bv_i^{(j)} \in A_\bu^{(j)}}
Z_i^{(j)} \mbox{ for } j=1,2, \mbox{ where}\cr
A_\bu^{(1)} = & \bigcup_{\bt \in \partial
B: \bn_\bt \cdot \bu > 0} \bt \times [0,\bn_\bt \cdot \bu)
\mbox{ and } A_\bu^{(2)} = \bigcup_{\bt \in \partial
B: \bn_\bt \cdot \bu < 0} \bt \times [\bn_\bt \cdot \bu,0).
\end{aligned} \leqno(2.3)$$
We define
$$Y(\bu) = Y^{(1)}(\bu)- Y^{(2)}(\bu) \mbox{ for
all } \bu \in {\bf R}^d. \leqno(2.4)
$$

\bigskip {\bf 3. First expression of asymptotic tail probability}. A key idea here is
a change of measure argument that allows us to obtain,
in Lemma 1, the tail probability of the maxima over
a local domain. To obtain the global probabilities in Theorem 1
from these local probabilities, we adapt the
Pickands-Quall-Watanabe technique  
from the Gaussian random field literature. Hence the
characterization of the constant $K$ in Theorem 1 bears a striking
resemblence to constants seen in the earlier
papers on Gaussian random fields 
though the distribution of $Y(\bu)$ here is
compound Poisson rather than Gaussian. 

Let $Q_\lambda$ be a probability measure under which
$\bX$ is a nonhomogeneous marked Poisson process with
rate $\lambda M(\theta_c)$ and mark distribution 
$F_c$ inside $B$, and rate $\lambda$ and mark distribution
$F$ outside $B$. Hence under $Q_\lambda$, for any set $A$,
Borel subset of $\mR^{d+1}$, $\# \{ i: (\bt_i,X_i) \in A \}$
follows a Poisson distribution with mean $\lambda M(\theta_c)
\int_{A \cap (B \times \mR)} d \bt \times dF_c(x) +
\lambda \int_{A \cap (B^c \times \mR)} d \bt \times dF(x)$
while $\# \{ i: (\bt_i,X_i) \in A \}$ and $\# \{ i: (\bt_i,
X_i) \in C \}$ are independent random variables for disjoint
sets $A$ and $C$. By (2.2),
$${dQ_\lambda \over dP_\lambda}(\bX) = \exp \{
\theta_c S(B)-\lambda \sigma_d(B)[M(\theta_c)-1] \}.
\leqno(3.1)
$$
Let $E_{\bt,m,\lambda}= \{ \sup_{\bt \prec \bv \prec \bt+
m\lambda^{-1} {\bf 1}}S(\bv+B) \geq \lambda c \}$. 
In the proof of Lemma 1 below, we analyze the event 
$E_{{\bf 0},m,\lambda}$ under $Q_\lambda$ before
applying the identity $P_\lambda(E_{{\bf 0},
m,\lambda}) = E_{Q_\lambda}[(dP_\lambda/
dQ_\lambda){\bf 1}_{E_{{\bf 0},m,\lambda}}]$. We shall
now define some terms required for the statement
of Lemma 1. 

For given $c > \max \{ 0, \mu \sigma_d(B) \}$, let 
$$I (=I_c)= \theta_c c - \sigma_d(B)[M(\theta_c)-1]. 
\leqno(3.2)
$$
It follows from Theorem 1 below that $I= - 
\lim_{\lambda \rightarrow \infty} \lambda^{-1} 
\log p_\lambda$ and hence $I$ is the large deviation
rate of the tail probability. If there exists $\eta >0$
such that $F$ is concentrated on $\pm \eta, \pm 2 \eta, \ldots$,
then we say that $F$ is arithmetic.
The largest $\eta$ with this property will be called
the span of $F$,
cf. Feller (1971) Section 5.2. If such $\eta$ does not exists, then we say
that $F$ is nonarithmetic. Let $\lfloor \cdot \rfloor$ denote the
greatest integer function and $''$ the second derivative of
a function.

\medskip {\sc Lemma} 1. {\it Let $c > \max \{ 0, \mu \sigma_d(B) \}$.
Define $x_\lambda=\theta_c(\lambda c-\eta \lfloor \lambda
c/\eta \rfloor)$ if $F$ is arithmetic with span $\eta$
and $x_\lambda=0$ if $F$ is nonarithmetic. Then for
all $\bt \in D$,
$$P_\lambda(E_{\bt,m,\lambda}) = 
P_\lambda(E_{{\bf 0},m,\lambda}) \sim [2 \pi \lambda \sigma_d(B)
M''(\theta_c)]^{-1/2} e^{-\lambda I+x_\lambda}
K_m \mbox{ as } \lambda \rightarrow \infty, \leqno(3.3)
$$
where
$$K_m = \begin{cases} \eta \Big[ (1-e^{-\eta \theta_c})^{-1}+
\sum_{\ell \in \eta {\mZ}^+} e^{\theta_c \ell} P \{
\sup_{{\bf 0} \prec \bu \prec m {\bf 1}} Y(\bu) \geq
\ell \} \Big] & \mbox{ if } F \mbox{ is arithmetic with span } \eta, \cr
\theta_c^{-1}+ \int_0^\infty e^{\theta_c y} P \{  
\sup_{{\bf 0} \prec \bu \prec m {\bf 1}} Y(\bu) \geq
y \} dy & \mbox{ if } F \mbox{ is nonarithmetic}. \end{cases}
\leqno(3.4)
$$
}

{\sc Proof}. By stationarity, $P_\lambda
(E_{\bt,m,\lambda})=P_\lambda(E_{{\bf 0},m,\lambda})$.
Let us first consider the case $F$
arithmetic with span 1. Then 
$$\begin{aligned} & P_\lambda(E_{{\bf 0},m,\lambda}) = 
P_\lambda \{ \sup_{{\bf 0} \prec \bv \prec m \lambda^{-1}
{\bf 1}} S(\bv+B) \geq \lambda c \} \cr
= & P_\lambda \{ S(B) \geq 
\lfloor \lambda c \rfloor \} 
+ \sum_{\ell=1}^\infty P_\lambda \{ S(B)=\lfloor \lambda 
c \rfloor - \ell , \sup_{{\bf 0} \prec \bv \prec 
m \lambda^{-1} {\bf 1}}  
[S(\bv+B)-S(B)] \geq \ell \}.
\end{aligned} \leqno(3.5)
$$
Let $B_a = \bigcap_{{\bf 0} \prec \bv \prec a {\bf 1}}
(\bv+B)$. Since $S(B \setminus B_{m \lambda^{-1}})$
and $\sup_{{\bf 0} \prec \bv \prec m \lambda^{-1}
{\bf 1}} [S(\bv+B)-S(B)]$ are functions of the
marked Poisson process occurring outside 
$B_{m \lambda^{-1}}$ and hence independent of
$S(B_{m \lambda^{-1}})$ under $Q_\lambda$, it
follows from (3.1) that 
$$\begin{aligned} & P_\lambda \{ S(B)=\lfloor \lambda c
\rfloor-\ell, \sup_{{\bf 0} \prec \bv \prec m \lambda^{-1}
{\bf 1}} [S(\bv+B)-S(B)] \geq \ell \} \cr
& \quad = e^{-\lambda I+x_\lambda+\theta_c \ell} Q_\lambda 
\{ S(B) = \lfloor \lambda c \rfloor-\ell, 
\sup_{{\bf 0} \prec \bv \prec m \lambda^{-1} {\bf 1}} 
[S(\bv+B)-S(B)] \geq \ell \} \cr
& \quad = e^{-\lambda I+x_\lambda+\theta_c \ell} 
\sum_{k=0}^\infty Q_\lambda 
\{ S(B_{m \lambda^{-1}}) = \lfloor \lambda c \rfloor-\ell
-k \} \cr
& \qquad \times Q_\lambda \{ S(B \setminus B_{m \lambda^{-1}})=k, 
\sup_{{\bf 0} \prec \bv \prec m \lambda^{-1} {\bf 1}} 
[S(\bv+B)-S(B)] \geq \ell \}. \end{aligned} \leqno(3.6)
$$
It follows from the local central limit theorem 
that for each $\ell \in {\mZ}$,
$$\begin{aligned} Q_\lambda \{ S(B_{m \lambda^{-1}}) = \lfloor \lambda
c \rfloor-\ell-k \} \leq & [1+o(1)] 
Q_\lambda \{ S(B)=\lfloor \lambda c \rfloor
-\ell \} \cr
\sim & [2 \pi \lambda \sigma_d(B) M''(\theta_c)]^{-1/2}
\mbox{ as } \lambda \rightarrow \infty, \end{aligned}
\leqno(3.7)$$
uniformly over $k \geq 0$, with $\leq$ replaced by
$=$ if we look at (3.7) with $k$ fixed. Hence by
(3.6) and (3.7),
$$\begin{aligned} & P_\lambda \{ S(B)=\lfloor \lambda c
\rfloor-\ell, \sup_{{\bf 0} \prec \bv \prec m \lambda^{-1}
{\bf 1}} [S(\bv+B)-S(B)] \geq \ell \} \cr
& \qquad \sim
[2 \pi \lambda \sigma_d(B) M''(\theta_c)]^{-1/2} e^{-\lambda I
+x_\lambda+\theta_c \ell} Q_\lambda \{
\sup_{{\bf 0} \prec \bv \prec m \lambda^{-1} {\bf 1}}
[S(\bv+B)-S(B)] \geq \ell \}. \end{aligned} \leqno(3.8)
$$
By (3.8) and the weak convergence of $\{ S(\lambda^{-1} \bu +B)-S(B):  
{\bf 0} \prec \bu \prec m {\bf 1} \}$ 
under $Q_\lambda$ to $\{ Y(\bu): {\bf 0} \prec \bu
\prec m {\bf 1} \}$ as 
$\lambda \rightarrow \infty$; see (2.3) and (2.4),   
$$\begin{aligned} & \sum_{\ell=1}^\infty P_\lambda \{
S(B) = \lfloor \lambda c \rfloor-\ell, \sup_{{\bf 0} \prec
\bv \prec m \lambda^{-1} {\bf 1}}[S(\bv+B)-S(B)] \geq
\ell \} \cr
& \qquad \sim e^{-\lambda I+x_\lambda} [2 \pi \lambda 
\sigma_d(B) M''(\theta_c)]^{-1/2}
\sum_{\ell=1}^\infty e^{\theta_c \ell} P
\{ \sup_{{\bf 0} \prec \bu \prec m {\bf 1}} Y(\bu) \geq
\ell \}. \end{aligned} \leqno(3.9)
$$
By a similar application of (3.1) and (3.7),
$$P_\lambda \{ S(B) \geq \lfloor \lambda c \rfloor \} \sim
[2 \pi \lambda \sigma_d(B) M''(\theta_c)]^{-1/2} e^{-\lambda I
+x_\lambda} \sum_{\ell=-\infty}^0 e^{\theta_c \ell}. \leqno(3.10)
$$
Substitution of (3.9) and (3.10) into (3.5) then proves
Lemma 1 when $F$ is arithmetic with span 1. For
$F$ arithmetic with arbitrary span $\eta$, we prove
Lemma 1 by replacing the sums in (3.5), (3.6), (3.9) and 
(3.10) by $\sum_{\ell \in \eta {\mZ}^+}$,
$\sum_{k \geq 0, k \in \eta \mZ}$ or 
$\sum_{\ell \leq 0, \ell \in \eta {\mZ}}$. For
nonarithmetic $F$, the sums are replaced by 
corresponding integrals. The detailed arguments 
are similar to the proof above and shall be omitted.
$\Box$

\medskip {\sc Theorem} 1. {\it Let $c > \max \{ 0, \mu \sigma_d(B) \}$ 
and define $x_\lambda$ as in Lemma 1. Then 
$$K := \lim_{m \rightarrow \infty} m^{-d} K_m \mbox{\it \ is
a well-defined positive and finite constant}. \leqno(3.11)
$$
Moreover, 
$$p_\lambda = P_\lambda \{ \sup_{\bv \in D} S(\bv+B) \geq 
\lambda c \} \sim [2 \pi \sigma_d(B) M''(\theta_c)]^{-1/2}
e^{-\lambda I+x_\lambda} \lambda^{d-1/2} \sigma_d(D) K \mbox{ as }
\lambda \rightarrow \infty. \leqno(3.12)
$$
}

{\sc Remarks}. By Jordan measurability of $D$, 
$$\# \Big\{ {\bf k} \in (a \bZ)^d: \prod_{j=1}^d
[k_j,k_j+a) \subset D \Big\} \sim \# \Big\{ {\bf k} \in (a \bZ)^d: 
\prod_{j=1}^d [k_j,k_j+a) \cap D \neq \emptyset \Big\} 
\mbox{ as } a \rightarrow 0. \leqno(3.13)
$$
The relation (3.12) still holds if $D$ is replaced by domains $D_\lambda$
that depends on $\lambda$, provided (3.13) holds with $D$
replaced by $D_\lambda$ and $a$ replaced by $m \lambda^{-1}$,
with limit $\lambda \rightarrow \infty$ for all large
$m$, and
$$\lim_{\lambda \rightarrow \infty} \lambda^{-1}
\log[\sigma_d(D_\lambda)] = 0. \leqno(3.14)
$$
Without condition (3.14), the correct relation is
$$P_\lambda \{ \sup_{\bv \in D_\lambda} S(\bv+B) \geq
\lambda c \} \sim 1 -\exp \{ - 
[2 \pi \sigma_d(B) M''(\theta_c)]^{-1/2} e^{-\lambda I+x_\lambda}
\lambda^{d-1/2} \sigma_d(D_\lambda) K \}. 
$$

\medskip We will now discuss an interesting case
of Theorem 1. In Example 1, we consider 
rectangular scanning sets on a marked Poisson random
field. We show here that an overshoot constant
derived from $F$ plays an important role in the tail
approximations. When $F$ is degenerate at 1, that is for
Poisson random fields rather than marked Poisson random
fields, the overshoot constant is equal to 1 and
disappears from the resulting formula.  

\medskip {\sc Example} 1. Let $B = \prod_{k=1}^d [0,b_k]$ with
$b_k > 0$ for all $k$. Since $\partial B$ is a union of $2d$ 
faces, with a pair of them orthogonal to each co-ordinate
vector, by (2.3) and (2.4),
$Y(\bu) = \sum_{k=1}^d [Y_k^{(1)}(u_k)
- Y_k^{(2)}(u_k)]$, where $Y_1^{(1)},\ldots,
Y_d^{(1)}, Y_1^{(2)},\ldots,Y_d^{(2)}$ are 
independent one-dimensional compound point processes.
The process $Y_k^{(1)}$, $1 \leq k \leq d$, 
is constructed from a marked
Poisson process having Poisson
rate $\prod_{\ell \neq k} b_\ell$; the surface area of the
face of $B$ orthogonal to the $k$th co-ordinate vector,
and mark distribution $F$. The process $Y_k^{(2)}$, 
$1 \leq k \leq d$, is
constructed from a marked Poisson process with 
Poisson rate $M(\theta_c) \prod_{\ell \neq k} b_\ell$
and mark distribution $F_c$. If $X$ is a random variable
with distribution $F_c$, we shall let $\bar F_c$ denote the
distribution of $-X$. We define $\bar F$ in a 
similar manner. Consider first $F$ nonarithmetic
and let $Y_k = Y_k^{(1)}-Y_k^{(2)}$. Then by (3.4)
and $P \{ \sup_{{\bf 0} \prec \bu \prec m {\bf 1}}
[\sum_{k=1}^d Y_k(u_k)] \geq y \}=1$ for $y \leq 0$,
$$\begin{aligned} & K_m =  
\int_{-\infty}^\infty e^{\theta_c y} P \Big\{ 
\sup_{{\bf 0} \prec \bu \prec m {\bf 1}} \Big[
\sum_{k=1}^d Y_k(u_k) \Big] \geq y \Big\} dy \cr
= & \int_{-\infty}^\infty \Big( \int_{-\infty}^y 
e^{\theta_c u} du \Big)
P \Big\{ \sup_{{\bf 0} \prec \bu \prec m {\bf 1}} 
\Big[ \sum_{k=1}^d Y_k(u_k) \Big] \in 
dy \Big\} = 
\theta_c^{-1} E \exp \Big[ \theta_c \sum_{k=1}^d
\sup_{0 \leq u_k \leq m} Y_k(u_k) \Big] \cr
= & \theta_c^{-1} \prod_{k=1}^d E \exp [ 
\theta_c \sup_{0 \leq u_k \leq m} Y_k(u_k) ] 
= \theta_c^{d-1} \prod_{k=1}^d \int_{-\infty}^\infty 
e^{\theta_c y} P \{
\sup_{0 \leq u_k \leq m} Y_k(u_k) \geq y \} dy. \end{aligned}
\leqno(3.15)
$$
Since $Y_k^{(1)}$ and $Y_k^{(2)}$ are independent
compound Poisson processes, it follows that 
$Y_k(u_k) = \sum_{j=1}^{N_k(u_k)} U_{kj}$, where
$N_k$ is a Poisson process with rate
$(\prod_{\ell \neq k} b_\ell)[1+M(\theta_c)]$ and
$U_{k1}, U_{k2}, \ldots$ are i.i.d. random variables
independent of $N_k$ such that 
$$P \{ U_{k1} \in du \} = [M(\theta_c)\bar F_c(du)+
F(du)]/[1+M(\theta_c)]. \leqno(3.16)
$$
Let $P_*$ be a probability measure under which
the distribution of $N_k$ is unchanged and
$U_{k1}, U_{k2},\ldots$ are i.i.d. random variables 
independent of $N_k$ satisfying
$$P_* \{ U_{k1} \in du \} = [\bar F(du) + M(\theta_c)
F_c(du)]/[1+M(\theta_c)]. \leqno(3.17)
$$
By (2.2), (3.16) and (3.17), 
$$(dP_*/dP)(U_{k1}) = e^{\theta_c U_{k1}}. \leqno(3.18)
$$
Suppressing the notation $k$, 
let $R_\ell =U_1+\cdots+U_\ell$ and 
$\tau_y = \inf \{ \ell \geq 1:
R_\ell \geq y \}$. Define the overshoot constant
$$\nu_c = \lim_{y \rightarrow \infty} E_* 
e^{-\theta_c(R_{\tau_y}-y)}, \leqno(3.19)
$$
where $E_*$ denotes expectation with respect to $P_*$.
See Siegmund (1985) Chapter 8 for the existence and
computation of $\nu_c$. By (3.17)-(3.19),
$$\begin{aligned} & \int_{-\infty}^\infty e^{\theta_c y} P \{
\sup_{0 \leq u_k \leq m} Y_k(u_k) \geq y \} dy = 
\theta_c^{-1} + \int_0^\infty
E_*[e^{\theta_c(y-R_{\tau_y})} {\bf 1}_{\{ 
\sup_{0 \leq u_k \leq m} Y_k(u_k) \geq y \}}] dy \cr
& \quad \sim \nu_c E_*[\sup_{0 \leq u_k \leq m} Y_k(u_k)]
\sim \nu_c m[c \sigma_d^{-1}(B)- \mu] \prod_{\ell \neq k} b_\ell,
\end{aligned} \leqno(3.20)
$$
noting that by (2.2), under $F_c$, $E_c X_1 = M'(\theta_c)/
M(\theta_c)=c/[\sigma_d(B) M(\theta_c)]$ and by 
definition, under $F$, $EX_1=\mu$. 
Substituting (3.20) into (3.15) and (3.11),
(3.12) then gives us  
$$p_\lambda \sim [2 \pi \sigma_d(B) M''(\theta_c)]^{-1/2}
e^{-\lambda I+x_\lambda} \lambda^{d-1/2}
\sigma_d(D) \{ \nu_c [c \sigma_d^{-1}(B)-\mu] \}^d
\Big( \chi_c \prod_{k=1}^d b_k \Big)^{d-1}, 
\leqno(3.21)
$$
where $\chi_c=\theta_c$ when $F$ is nonarithmetic. 
Using similar arguments, 
the relation (3.21) can also be shown to hold for $F$
arithmetic with span $\eta$, by defining $\nu_c$ in
(3.19) with limit $y \in \eta {\mZ}$, $y \rightarrow 
\infty$ and $\chi_c = \eta^{-1}(1-e^{-\eta \theta_c})$. 

\bigskip
{\bf 4. An alternative approach}. The evaluation
of the constant $K$ in Example 1 for rectangular kernels follows along the lines of
Hogan and Siegmund (1986). However, when $B$ is not rectangular,
the expression of $K$ via (3.4) and (3.11) does not seem to be helpful
except for indicating how the proofs of Lemma 1 and Theorem 1 is expected to
proceed. This is unsatisfying since kernels
of other shapes are often used in practice. For example, in epidemiology
and geostatistical applications, the circular kernel $B = \{ (t_1,t_2):
t_1^2+t_2^2 \leq 1 \}$ provides a more desirable co-ordinate free
space symmetry. For space-time problems, the corresponding kernel is 
the cylindrical scanning set $B = \{ (t_1,t_2,t_3):
t_1^2+t_2^2 \leq 1, |t_3| \leq 1 \}$. 

To search for an alternative formulation of $K$, it is
best to start with the special case $F$ concentrated at 1, for which the identity
$$K = E[\sigma_d^{-1} (\{ \bu \in {\bf R}^d: Y(\bu)=0 \})] \leqno(4.1)
$$
holds. This identity looks surprising initially 
because the right hand side involves only the occupation
measure of the conditional process $Y$ at 0 and does not seem to be related
to the maxima of $Y$. Is (4.1) true for general $F$ ? Before answering this
question, we first show how (4.1) can be utilized to provide a lower bound
for $K$. 

\medskip
\begin{table}
\begin{tabular}{c|cccccccccc}
$c$ & 2 & 3 & 4 & 5 & 6 & 7 & 8 & 9 & 10 & $\infty$ \cr \hline
(I) $F$ concentrated at 1 & .160 & .497 & .818 & 1.08 & 1.29 & 1.47 & 1.62 &
1.75 & 1.85 & 3.12 \cr \hline
Lower bound of (I) & .0235 & .0795 & .137 & .188 & .232 & .268 & .299 & 
.325 & .348 & .636 \cr \hline
(II) $F \sim N(0,1)$ & .153 & .324 & .495 & .654 & .797 & .929 & 1.04 &
1.15 & 1.25 & 3.12 
\end{tabular}
\caption{ \ Entries of $K/(1+c)^2$ for the second row and $K \theta_c/
[1+M(\theta_c)]^2$ for the fourth row for the kernel $B = \{ \bt: \| \bt \|
\leq 1 \}$ with $d=2$. These numbers have an approximate
1\% numerical error. The third row is obtained from the inequality $K \geq 2(c-1)^3/[\pi(1+c)]$,
see Example 2.}
\end{table}

\medskip {\sc Example 2}. Let $F$ be concentrated at 1. By (4.1),
$K \geq \{ E[\sigma_d(\{ \bu: Y(\bu)=0 \})] \}^{-1}$. 
Let $B_d = \{ \bt: \| \bt \| \leq 1 \}$ with $\bt \in {\bf R}^d$ and let $C_d = 
\sigma_{d-1}(\partial B_d)/\sigma^d_{d-1}(B_{d-1}) = [d \pi^{d/2}/\Gamma(d/2+1)]/
[ \pi^{(d-1)/2}/\Gamma((d+1)/2)]^d$. Then
$$\begin{aligned}
& E[\sigma_d(\{ \bu: Y(\bu)=0 \})]=C_d \int_0^\infty r^{d-1} \sum_{k=0}^\infty
e^{-r(1+c)}[c^k r^{2k}/(k!)^2] \ dr \cr
& \quad = C_d \sum_{k=0}^\infty c^k \Gamma(2k+d)/[(k!)^2(1+c)^{2k+d}] = 
[C_d/(1+c)^d] f_d^{(d-1)}(\sqrt{c}/(1+c)),
\end{aligned} \leqno(4.2)
$$
where $f_\ell(x)=x^\ell/\sqrt{1-4x^2}$ and $g^{(\ell)}$ is the $\ell$th
derivative of a function $g$. Similar computations can also be carried out for
kernels of other shapes.

\medskip
{\sc Theorem 2}. {\it Let $c > \max \{ 0, \mu \sigma_d(B) \}$. Then
$$K = \chi_c^{-1} E \{ [(1-\exp(\theta_c \sup \{ \bu \in {\bf R}^d:
Y(\bu) < 0 \})] \sigma_d^{-1}(\{ \bu \in {\bf R}^d: Y(\bu)=0 \}) \} \leqno(4.3)
$$ 
where $\chi_c = \eta^{-1}(1-e^{-\eta \theta_c})$ if $F$ is
arithmetic with span $\eta$ and $\chi_c=\theta_c$ if
$F$ is nonarithmetic. 
}

\medskip
Let $\{ (\bt(i), y_i): i \geq 1 \}$ be a unit rate Poisson process
 defined on $\partial B
\times [0,\infty)$ and define the random set
$$\Omega = \{ \bu \in {\bf R}^d: \bn_{\bt(i)} \cdot \bu < 0
\mbox{ or } y_i(\bn_{\bt(i)} \cdot \bu) \geq \| \bu \|^2
\mbox{ for all } i \}. 
$$
Then $[1+M(\theta)]^d \sigma_d(\{ \bu: Y(\bu)=0 \}) \Rightarrow \sigma_d(\Omega)$
as $c \rightarrow \infty$. Hence we obtain the following.

\medskip {\sc Corollary} 1. {\it If $F$ is nondegenerate, then
$K/\{ \chi_c^{-1} [1+M(\theta_c)]^d \}$ is bounded above by $E[\sigma_d^{-1}(\Omega)]$
for all $c$ and tends towards $E[\sigma_d^{-1}(\Omega)]$ as $c \rightarrow
\infty$. If $F$ is degenerate at $\eta > 0$, then $K/[1+M(\theta_c)]^d$ is bounded above by
$\eta E[\sigma_d^{-1}(\Omega)]$ for all $c$ and tends toward $\eta E[\sigma_d^{-1}(\Omega)]$
as $c \rightarrow \infty$. }

\medskip The case $F$ degenerate at $\eta$ stands out because $\sup \{ Y(\bu):
Y(\bu) < 0 \} = -\eta$ with probability 1. Note also that $\Omega$
depends only on the kernel $B$ and not $F$. 

\medskip
{\sc Proof of Theorem 2}. Let us first consider $F$ arithmetic
with span 1.  To simplify notations, select $\lambda$ such that $x_\lambda=0$
(i.e. $\lambda \in {\bf Z}/c$). We shall also abuse notation here 
and write $S(\bv)$ in place  of $S(\bv+B)$. By
the change of measure argument in the proof of Lemma 1 and 
the probability bounds obtained in Lemmas A.1 and A.2,
for any integers $0 \leq k < \ell$, $\bt$ in the interior of $D$ 
and $dw \in (0,\infty)$, 
$$\begin{aligned}
& P_\lambda \Big\{ \sup_{\bv \in D} S(\bv) = \lambda c, S(\bt)=
\lambda c-k, \sigma_d(\{ \bv: S(\bv) = \lambda c-k \}) \in
\lambda^{-d} dw, \cr
& \quad \sup \{ S(\bv): S(\bv) < S(\bt) \} = \lambda c-\ell \Big\} 
\sim [2 \pi \lambda \sigma_d(B) M''(\theta_c)]^{-1/2}
e^{-\lambda I+\theta_c k} \sigma_d(D) \cr
& \quad \times P \Big\{ \sup_{\bu \in
{\bf R}^d} Y(\bu) =k, \sigma_d(\{ \bu:Y(\bu)=0 \})
\in dw, \sup_{\bu \in {\bf R}^d} \{
Y(\bu): Y(\bu) < 0 \}=k-\ell \Big\}. \end{aligned} \leqno(4.4)
$$
Multiplying (4.4) by $\lambda^d (e^{-\theta_c k}-e^{-\theta_c \ell})/w$
and integrating over $\bt \in D$ and $w>0$, we obtain
$$\begin{aligned}
& (e^{-\theta_c k}-e^{-\theta_c \ell}) P_\lambda \Big\{
\sup_{\bv \in D} S(\bv)=\lambda c, \sigma_d(\{ \bv: S(\bv)=
\lambda c-k \})>0, \cr
& \quad \sigma_d(\{ \bv: S(\bv)=\lambda c-\ell \})
> 0, \sigma_d(\{ \bv: \lambda c -\ell < S(\bv) < \lambda c-k \})
=0 \Big\} \cr
\sim & [2 \pi \sigma_d(B) M''(\theta_c)]^{-1/2} e^{-\lambda I}
\lambda^{d-1/2} \sigma_d(D) (1-e^{\theta_c(k-\ell)}) \cr
& \quad \times 
E \Big[ \sigma_d^{-1}(\{ \bu: Y(\bu)=0 \}) {\bf 1}_{\{
\sup_{\bu \in {\bf R}^d} Y(\bu)=k, \sup_{\bu \in {\bf R}^d}
\{ Y(\bu): Y(\bu) < 0 \}=k-\ell \}} \Big].
\end{aligned} \leqno(4.5)
$$
Theorem 2 then follows by adding up (4.5) over the integers $0 \leq k < \ell$
and comparing against 
$$P_\lambda \Big\{ \sup_{\bv \in D} S(\bv+B) = \lambda c \Big\} \sim
[2 \pi \sigma_d(B) M''(\theta_c)]^{-1/2} (1-e^{-\theta_c})
e^{-\lambda I} \lambda^{d-1/2} \sigma_d(D) K, \leqno(4.6)
$$
a straightforward modification of Theorem 1. For $F$
arithmetic with arbitrary span $\eta > 0$ or $F$ nonarithmetic,
the arguments are similar. $\Box$ 

\medskip
{\bf 5. A relook at the Poisson clumping heuristic}. 
In this section, we consider a continuous valued random process $X(\bt)$,
$t \in D$. To make the discussion concrete, we pick the isotropic
mean zero Gaussian random field $X(\bt)$, $\bt \in {\bf R}^d$ satisfying
$$E[X(\bt) X(\bt+\bs)] \sim 1-a \| \bs \|^\alpha \mbox{ as } \| \bs \| \rightarrow
0 \leqno(5.1)
$$
for some $0 < \alpha \leq 2$ and $0 < a < \infty$. It was shown in 
Bickel and Rosenblatt (1973) and Qualls and Watanabe
(1973) that
$$P \{ \sup_{\bt \in D} X(\bt) \geq c \} \sim (2 \pi)^{-1/2} c^{2d/\alpha-1}
e^{-c^2/2} a^{d/\alpha} \wtd K \mbox{ as } c \rightarrow \infty. \leqno(5.2)
$$
The approach is via a conditioning on $X(\bt) = c-y/c$ for $y > 0$ which 
leads to the expression 
$$\wtd K = \lim_{m \rightarrow \infty} m^{-d} \int_0^\infty e^y P \{ \sup_{\bu
\in [0,m]^d} Y(\bu) \geq y \} \ dy, \leqno(5.3)
$$
where $Y$ is a Gaussian process satisfying 
$$E Y(\bu) = -\| \bu \|^\alpha, \quad \cov(Y(\bu),Y(\bv)) = \| \bu \|^\alpha
+\| \bv \|^\alpha-\| \bu-\bv \|^\alpha. \leqno(5.4)
$$
Aldous (1989) using the Poisson clumping heuristic,
conditioned instead on $X(\bu) \geq c+y/c$ for $y>0$ and it follows from
this approach that
$$\wtd K = E[\sigma_d^{-1}(\{ Y(\bu): Y(\bu) \geq -Z \})], \leqno(5.5)
$$
where $Z$ is an independent exponential random variable with mean 1. In
Theorem 3, we apply the technique used to prove Theorem 2 to provide a
differential form of (5.5). 

\medskip {\sc Theorem 3}. $\wtd K = \lim_{\xi \rightarrow 0} \int_0^\xi E[\sigma_d^{-1}(\{ \bu:
-b < Y(\bu) \leq \xi-b \})] db$.

\medskip {\sc Example 3}. We shall provide lower bounds of 
$\wtd K$ using the harmonic mean inequality
as in J20 of Aldous (1989). Let $B = \{ \bt: \| \bt \| \leq 1 \}$. By Theorem 3, 
$\wtd K \geq \lim_{\xi \rightarrow 0} \int_0^\xi \{ E[\sigma_d(\{ \bu:
-b < Y(\bu) \leq \xi-b \})] \}^{-1} \ db$
and hence
$$\wtd K^{-1} \leq \sigma_{d-1}(\partial B) \int_0^\infty
r^{d-1} (4 \pi r^\alpha)^{-1/2} e^{-r^{2 \alpha}/(4r^\alpha)} \ dr.
$$
This leads to the inequality
$$\wtd K \geq d^{-1} \pi^{(1-d)/2} 4^{1-d/\alpha} \alpha \Gamma \Big( \frac{d}{2}+1 \Big) \Big/
\Gamma \Big( \frac{d}{\alpha}
-\frac{1}{2} \Big). \leqno(5.6)
$$
In the case $\alpha=2$, $Y$ has a simple characterization from which 
$\wtd K=\pi^{-d/2}$ can be computed. For $d=2$, the right hand side of 
(5.6) is $\pi^{-1}$(=$\wtd K$)
and for $d=3$, it is $1/(4 \sqrt{\pi})$. 

\medskip {\sc Proof of Theorem 3}. 
Let $\xi > 0$, $0 \leq v < \xi$, $X_{\sup}= \sup_{\bt \in D} X(\bt)$ and
$Y_{\sup}=\sup_{\bu \in {\bf R}^d} Y(\bu)$. For any integer $k \geq 0$,
$\bt$ in the interior of $D$ and $dw \in (0,\infty)$,
$$\begin{aligned}
& P \{ X_{\sup} \geq c, X_{\sup}-[(k+1) \xi-v]/c < X(\bt) \leq X_{\sup}-
(k \xi-v)/c, \cr
& \quad \sigma_d(\{ \bu: X_{\sup}-[(k+1) \xi-v]/c < X(\bu) \leq
X_{\sup}-(k \xi-v)/c \}) \in (c^2 a)^{-d/\alpha} dw \} \cr
& \quad \sim (2 \pi)^{-1/2} c^{-1} e^{-c^2/2} \Big\{ \Big( \int_{-(k \xi-v)}^\infty
e^{-y} \ dy \Big) P \{ k \xi-v \leq Y_{\sup} < (k+1) \xi-v, \cr
& \quad \sigma_d(\{ \bu:
Y_{\sup}-[(k+1) \xi-v] < Y(\bu) \leq Y_{\sup}-(k \xi-v) \}) \in dw \} \cr
& \quad + \int_{-[(k+1) \xi-v]}^{-(k \xi-v)} e^{-y} P \{ y \leq Y_{\sup}
< (k+1) \xi-v, \cr
& \quad \sigma_d(\{ \bu: Y_{\sup}-[(k+1) \xi-c] < Y(\bu) \leq 
Y_{\sup}-(k \xi-c) \}) \in dw \} dy \Big\}. 
\end{aligned} \leqno(5.7)
$$
Multiply (5.7) by $(e^{-k \xi}-e^{-(k+1) \xi})/[(c^2 a)^{-d/\alpha} w]$,
then integrating over $\bt \in D$ and $w > 0$  and add over $k \geq 0$. Then  
$$\begin{aligned}
& \sum_{k=0}^\infty (e^{-k \xi}-e^{-(k+1) \xi}) P \{ X_{\sup} \geq c \} \cr
& \quad \sim (2 \pi)^{-1/2} c^{2d/\alpha-1} a^{d/\alpha} e^{-c^2/2} \Big\{ 
(e^{-v}-e^{-v-\xi}) \sum_{k=0}^\infty E[ {\bf 1}_{\{ k \xi-v \leq 
Y_{\sup} < (k+1) \xi-v \}} \cr
& \quad \times \sigma_d^{-1}(\{ \bu: Y_{\sup}-
[(k+1) \xi-v] < Y(\bu) \leq Y_{\sup}-(k \xi-v) \})] + o_\xi(1) \Big\}, \end{aligned}
$$
where $o_\xi(1) \rightarrow 0$ as $\xi \rightarrow 0$. By (5.2), 
$$
\wtd K = \lim_{\xi \rightarrow 0} \int_0^\xi E[\sigma_d^{-1}(\{ \bu: Y_{\sup}-
(\xi \lfloor (Y_{\sup}+v)/\xi \rfloor-v-\xi) < Y(\bu) \leq Y_{\sup}
-(\xi \lfloor (Y_{\sup}+v)/\xi \rfloor-v )\}) ] dv \leqno(5.8)
$$
and Theorem 3 is shown. $\Box$

\bigskip
\centerline{\sc Acknowledgements}

\smallskip I would like to thank an associate editor and a referee for their
valuable comments and reference. 

\bigskip 
\centerline{\sc Appendix: Proof of Theorem 1}

\bigskip Let
$J=[0,1]^d$, $\underline C_a = \{ \bk \in (a \bZ)^d:
\bk +aJ \subset D \}$ and $\overline C_a = \{ \bk \in
(a \bZ)^d: \bk+aJ \cap D \neq \emptyset \}$. Then
$\{ \bk+aJ: \bk \in \underline C_a \}$ and
$\{ \bk+aJ: \bk \in \overline C_a \}$ are lower and
upper coverings of $D$ respectively by cubes of length $a$.
We shall show via Lemmas A.1 and A.2 that  
$$\lim_{m \rightarrow \infty} \limsup_{\lambda \rightarrow \infty} 
\Big[ \sum_{\bu \in \underline C_{m \lambda^{-1}}} 
P_\lambda \Big\{ \bigcup_{\bw \in \underline C_{m \lambda^{-1}}:
\bw \neq \bu} (E_{\bu,m,\lambda} \cap E_{\bw,m,\lambda}) \Big\}/
(\lambda^{d-1/2} e^{-\lambda I}) \Big] = 0. \leqno({\rm A}.1)
$$
Let $f(\lambda)=[2 \pi \sigma_d(B) M''(\theta_c)]^{-1/2}
\lambda^{d-1/2} e^{-\lambda I+x_\lambda}$. Then by Lemma 1,
$$P_\lambda(E_{\bu,m,\lambda})
\sim K_m f(\lambda)/\lambda^d \mbox{ as } \lambda \rightarrow
\infty. \leqno({\rm A}.2)
$$ 

Given $\epsilon > 0$, let $m_\epsilon$ be large enough such
that for all $m \geq m_\epsilon$, the expression in the
square brackets on the left-hand
side of ({\rm A}.1) does not exceed $\epsilon$ for all large $\lambda$. Then
by ({\rm A}.2), for all $m \geq m_\epsilon$, 
$$\begin{aligned} & \liminf_{\lambda \rightarrow \infty} [
\lambda^{-d} (\# \underline C_{m/\lambda})
K_m - \epsilon \lambda^{d-1/2} e^{-\lambda I}/f(\lambda)] 
\leq \liminf_{\lambda \rightarrow
\infty} [p_\lambda/f(\lambda)] \cr
& \qquad \leq \limsup_{\lambda \rightarrow \infty}
[p_\lambda/f(\lambda)] \leq \limsup_{\lambda \rightarrow \infty} 
[\lambda^{-d} (\#
\overline C_{m/\lambda}) K_m]. \end{aligned} \leqno({\rm A}.3)
$$
Since $D$ is Jordan-measurable, 
$$\# \underline C_a \sim \# \overline C_a \sim a^{-d}
\sigma_d(D) \mbox{ as } a \rightarrow 0. \leqno({\rm A}.4)
$$
Noting that $\liminf_{\lambda \rightarrow \infty}
[p_\lambda/f(\lambda)]$ and $\limsup_{\lambda \rightarrow
\infty} [p_\lambda/f(\lambda)]$ are fixed real numbers and 
$x_\lambda$ is bounded, it follows from ({\rm A}.3) and ({\rm A}.4) 
that $m^{-d} K_m$ is
Cauchy. Hence $K=\lim_{m \rightarrow \infty} m^{-d}
K_m$ exists and (3.12) follows from ({\rm A}.3). 

We will now state and prove Lemmas A.1 and A.2 before
providing the complete proofs of both ({\rm A}.1) and Theorem 1. To avoid
repetitive arguments, we will state and prove all subsequent results
assuming $F$ is arithmetic with span 1. The modifications required
to extend these results to arbitrary $F$ are relatively 
straightforward and will not be discussed. 

\medskip {\sc Lemma} A.1. {\it 
$$\lim_{r \rightarrow \infty} \limsup_{\lambda \rightarrow \infty} 
\Big[ P_\lambda \{ S(B) < \lfloor \lambda c \rfloor -r, 
\sup_{{\bf 0} \prec \bv \prec \lambda^{-1} {\bf 1}} 
S(\bv+B) \geq \lfloor \lambda c \rfloor \}/(\lambda^{-1/2} e^{-\lambda I})
\Big] = 0. \leqno({\rm A}.5)
$$
}

{\sc Proof}. By (3.8) with $m=1$ and the weak convergence of
$S(\lambda^{-1} \bu+B)-S(B)$ to $Y(\bu)$ under $Q_\lambda$,
$$\begin{aligned} & P_\lambda \{ S(B) < \lfloor \lambda c \rfloor-r, 
\sup_{{\bf 0} \prec \bv \prec \lambda^{-1} {\bf 1}} S(\bv+B) \geq
\lfloor \lambda c \rfloor \} \cr
& \quad \sim [2 \pi \lambda \sigma_d(B) M''(\theta_c)]^{-1/2}
e^{-\lambda I+x_\lambda} \sum_{\ell =r+1}^\infty
e^{\theta_c \ell} P \{ \sup_{{\bf 0} \prec \bu \prec {\bf 1}}
Y(\bu) \geq \ell \} \mbox{ as } \lambda \rightarrow
\infty. \end{aligned} \leqno({\rm A}.6)
$$
Let $x^+=\max \{ x,0 \}$ and $x^- = \max \{ -x,0 \}$.
It follows from (2.3) and (2.4) that
$$\sup_{{\bf 0} \prec \bu \prec {\bf 1}} Y(\bu) \leq Z^* :=
\sum_{i:\bv_i^{(1)} \in A^{(1)}} [Z_i^{(1)}]^+ +
\sum_{i:\bv_i^{(2)} \in A^{(2)}} [Z_i^{(2)}]^-,
\leqno({\rm A}.7)
$$
where $A^{(1)} = \partial B \times [0,d^{1/2})$ and 
$A^{(2)} = \partial B \times [-d^{-1/2},0)$. We can also express
$Z^* = \sum_{j=1}^N V_j$, where $N$ is a Poisson random variable with
mean $\kappa := 1-F(0)+M(\theta_c)F_c(0)$
and $V_1, V_2, \ldots$ are i.i.d. random variables independent
of $N$ with $g:=E e^{\tilde \theta V_1} < \infty$ for
some $\wtd \theta > \theta_c$. Since 
$E e^{\tilde \theta Z^*} = e^{\kappa (g-1)}$, it follows 
from Chebyshev's inequality that
$P \{ Z^* \geq \ell \} \leq e^{\kappa(g-1)-\tilde
\theta \ell}$. Hence by ({\rm A}.7), 
$$\sum_{\ell=r+1}^\infty e^{\theta_c \ell} P_\lambda
\{ \sup_{{\bf 0} \prec \bu \prec {\bf 1}} Y(\bu) \geq \ell \} =
O(e^{\kappa(g-1)-(\tilde \theta-\theta_c)r}). \leqno({\rm A}.8)
$$
Lemma 2 follows from ({\rm A}.6), ({\rm A}.8) and because $x_\lambda$
is bounded. $\Box$

\medskip
{\sc Lemma} A.2. {\it Let $r \geq 0$ and $L > 0$ be given. Then
$$\lim_{k \rightarrow \infty} \limsup_{\lambda \rightarrow \infty} \Big[ 
\sum_{\bv \in (\lambda^{-1} \mZ)^d: k \lambda^{-1} \leq 
\| \bv \|_\infty \leq L} P_\lambda \{ S(B) \geq 
\lfloor \lambda c \rfloor -r, S(\bv+B) \geq \lfloor 
\lambda c \rfloor -r \}/(\lambda^{-1/2}
e^{-\lambda I}) \Big] = 0. \leqno({\rm A}.9)
$$
}

{\sc Proof}. Let $\wQ_\lambda(=\wQ_{\lambda,\bv})$ 
be the probability measure under which the
marked Poisson process $\bX$ has Poisson rate $\lambda M(\theta_c)$ on
$B_1 := B \cap (\bv+B)$, rate $\lambda M(\theta_c/2)$ on
$B_2:= (B \setminus (\bv+B)) \cup ((\bv+B) \setminus B)$ and
rate $\lambda$ elsewhere on ${\mR}^d$. Moreover 
we require that under $\wQ_\lambda$, 
the marks have distribution
$F_c$ on $B_1$, distribution $\wtd F$ satisfying $\wtd F(dx) = 
e^{\theta_c x/2} F(dx)/M(\theta_c/2)$ on $B_2$
and $F$ elsewhere on ${\mR}^d$. Then
$${d \wQ_\lambda \over dP_\lambda}
(\bX) = \Big( \prod_{i:\bt_i \in B_1} e^{\theta_c X_i}
\Big) e^{-\lambda \sigma_d(B_1)[M(\theta_c)-1]}
\Big( \prod_{i:\bt_i \in B_2} e^{\theta_c X_i/2} \Big) e^{-\lambda
\sigma_d(B_2)[M(\theta_c/2)-1]}. \leqno({\rm A}.10)
$$
Since $M$ is a convex function and  
$M(0)=1$, $\zeta:=[M(\theta_c)-1]-
2[M(\theta_c/2)-1] > 0$. We can thus express ({\rm A}.10) as
$${d \wQ_\lambda \over dP_\lambda}(\bX) = e^{(\theta_c/2)[S(B)+S(\bv+B)]-
\lambda[\sigma_d(B)+\sigma_d(\bv+B)][M(\theta_c)-1]/2+\lambda
\zeta \sigma_d(B_2)/2}, 
$$
and it follows from (3.2) and an analogue of (3.7) that 
$$\begin{aligned} & P_\lambda \{ S(B) \geq \lfloor 
\lambda c \rfloor-r, 
S(\bv+B) \geq \lfloor \lambda c \rfloor-r \} \leq  P_\lambda 
\{ S(B)+S(\bv+B) \geq 2(\lfloor \lambda c
\rfloor-r) \} \cr
& \quad =
E_{\tilde Q_\lambda} \Big[ {dP_\lambda \over d 
\wQ_\lambda} {\bf 1}_{\{ S(B)+S(\bv+B) \geq
2(\lfloor \lambda c \rfloor-r) \}} \Big] = 
O(e^{-\lambda I-\lambda \zeta \sigma_d(B \setminus
(\bv+B))} \lambda^{-1/2}). \end{aligned} \leqno({\rm A}.11)
$$

Let $\| \be \|=1$ and $\Pi_\be = \{ \bb-(\be \cdot \bb) \be: \bb \in B \}$ 
the projected surface of $B$ on a $(d-1)$-dimensional
hyperplane orthogonal to $\be$. Then  
$\beta := \inf_{\| \be \| = 1} \sigma_{d-1}(\Pi_\be) > 0$. 
Hence there exists $\epsilon > 0$ such that 
$$\sigma_d(B \setminus (\bv+B)) \geq \| \bv \| \beta/2
\geq \| \bv \|_\infty \beta/2 \mbox{ for all }
\| \bv \|_\infty \leq \epsilon. \leqno({\rm A}.12)
$$
By ({\rm A}.11) and ({\rm A}.12), it follows that  
$$\begin{aligned} & \sum_{\bv \in (\lambda^{-1} \mZ)^d:
k \lambda^{-1} \leq \| \bv \|_\infty \leq  \epsilon}
P_\lambda \{ S(B) \geq \lfloor \lambda c \rfloor-r, 
S(\bv+B) \geq \lfloor \lambda c \rfloor-r \} \cr
& \qquad = O \Big( \lambda^{-1/2} e^{-\lambda I}
\sum_{\ell \geq k} \ell^{d-1} e^{-\zeta \ell \beta/2}
\Big) = O(\lambda^{-1/2} e^{-\lambda I} k^{d-1} 
e^{-\zeta k \beta/2}).
\end{aligned} \leqno({\rm A}.13)
$$
Moreover, since $\alpha:= \inf_{\| \bv \|_\infty > \epsilon}
\sigma_d(B \setminus (\bv+B)) > 0$, it follows from 
({\rm A}.11) that
$$\sum_{\bv \in (\lambda^{-1} \mZ)^d: 
\epsilon < \| \bv \|_\infty \leq L} 
P_\lambda \{ S(B) \geq \lfloor \lambda c \rfloor-r, 
S(\bv+B) \geq \lfloor \lambda c \rfloor-r \} 
= O( \lambda^{d-1/2} e^{- \lambda I -
\lambda \zeta \alpha}). \leqno({\rm A}.14)
$$
Lemma 3 then follows from combining ({\rm A}.13) and ({\rm A}.14). $\Box$

\medskip {\sc Proof of} ({\rm A}.1). Let $\epsilon > 0$. By Lemma A.1
and stationarity, we can select $r$ large enough such that
$$\gamma_{\bu,\lambda} := P_\lambda \{ S(\bu+B) < \lfloor \lambda
c \rfloor-r, \sup_{\bu \prec \bv \prec \bu+{\bf 1}} S(\bv+B)
\geq \lfloor \lambda c \rfloor \} \leq \epsilon \lambda^{-1/2}
e^{-\lambda I} \leqno({\rm A}.15)
$$
for all large $\lambda$. Let $k=\lfloor m^{1/2} \rfloor$,
$\Gamma_m = \{ \bt \in \mZ^d: k {\bf 1} \prec \bt \prec
(m-k) {\bf 1} \}$ and $\Omega_m = \{ \bt \in \mZ^d: {\bf 0}
\prec \bt \prec m {\bf 1} \} \setminus \Gamma_m$. Then
$$\begin{aligned} & \sum_{\bu \in \underline C_{m \lambda^{-1}}}
P_\lambda \Big\{ \bigcup_{\bw \in \underline C_{m \lambda^{-1}}:
\bw \neq \bu}
(E_{\bu,m,\lambda} \cap E_{\bw,m,\lambda}) \Big\} \cr
& \qquad \leq
\sum_{\bu \in \underline C_{\lambda^{-1}}} \gamma_{\bu,\lambda}
+ \sum_{\bu,\bw \in \underline C_{m \lambda^{-1}}:\bw \neq \bu}
P_\lambda(G_{\bu,m,\lambda} \cap G_{\bw,m,\lambda}) + \sum_{\bu \in
\underline C_{m \lambda^{-1}}} P_\lambda(H_{\bu,m,\lambda}), 
\end{aligned} \leqno({\rm A}.16)
$$
where
$$G_{\bu,m,\lambda} = P_\lambda \{
S(\bv+B) \geq \lfloor \lambda c \rfloor-r \mbox{ for some }
\bv \in \bu + \lambda^{-1} \Gamma_m \}, \leqno({\rm A}.17)
$$
$$H_{\bu,m,\lambda} = P_\lambda \{ S(\bv+B) \geq \lfloor \lambda c
\rfloor-r \mbox{ for some } \bv \in \bu + \lambda^{-1} \Omega_m \}. 
\leqno({\rm A}.18)
$$
By ({\rm A}.4) and ({\rm A}.15),
$$\sum_{\bu \in \underline C_{\lambda^{-1}}} \gamma_{\bu,\lambda}
\leq [\epsilon+o(1)] \lambda^{d-1/2} e^{-\lambda I} \sigma_d(D).
\leqno({\rm A}.19)
$$
Let $L > \sup_{\bx,\by \in D} \| \bx - \by \|_\infty$. Then by
({\rm A}.17), stationarity and Lemma A.2, there exists $m$ large enough such
that for all $\bu \in \underline C_{m \lambda^{-1}}$ and
large $\lambda$,
$$\begin{aligned} & \sum_{\bw \in \underline C_{m \lambda^{-1}}:\bw \neq
\bu} P_\lambda
(G_{\bu,m,\lambda} \cap G_{\bw,m,\lambda}) \cr
& \leq 
m^d \sum_{\bv \in (\lambda^{-1} {\mZ})^d: k \lambda^{-1}
\leq \| \bv \|_\infty \leq L} P_\lambda \{ S(B) \geq \lfloor
\lambda c \rfloor-r, S(\bv+B) 
\geq \lfloor \lambda c \rfloor
-r \} 
\leq [\epsilon+o(1)] m^d \lambda^{-1/2} e^{-\lambda I}. \end{aligned}
$$
Hence by ({\rm A}.4), 
$$\sum_{\bu,\bw \in \underline C_{m \lambda^{-1}}: \bw
\neq \bu}
P_\lambda (G_{\bu,m,\lambda} \cap G_{\bw,m,\lambda})
\leq \sigma_d(D)[\epsilon+o(1)] \lambda^{d-1/2} e^{-\lambda I}.
\leqno({\rm A}.20)
$$
Since $k = \lfloor m^{1/2} \rfloor$, 
$\# \Omega_m = (m+1)^d -(m+1-2k)^d = O(m^{d-1/2})$ and it 
follows from ({\rm A}.4) and 
a modification of (3.10) (with $\sum_{\ell=-\infty}^r$
instead of $\sum_{\ell=-\infty}^0$) that
$$\sum_{\bu \in \underline C_{m \lambda^{-1}}} P_\lambda
(H_{\bu,m,\lambda}) = O(m^{d-1/2}(\lambda m^{-1})^d \sigma_d(D) \lambda^{-1/2}
e^{-\lambda I}) \leq \epsilon \lambda^{d-1/2} e^{-\lambda I} \sigma_d(D)
\leqno({\rm A}.21)
$$
for all large $m$. We then obtain
({\rm A}.1) from ({\rm A}.16) and ({\rm A}.19)-({\rm A}.21) by choosing $\epsilon$ 
arbitrarily small. $\Box$

\medskip {\sc Proof of Theorem 1}. By the arguments in the beginning
of Section 4, it remains to show that $K$ is positive and finite.
By Lemma A.1, there exists $r$ large enough such that
$$P_\lambda \{ S(B) < \lfloor \lambda c \rfloor -r, \sup_{{\bf 0}
\prec \bv \prec \lambda^{-1} {\bf 1}} S(\bv+B) \geq \lfloor
\lambda c \rfloor \} \leq \lambda^{-1/2} e^{-\lambda I} \leqno({\rm A}.22)
$$
for all large $\lambda$. 
Moreover, by a modification of (3.10) (with $\sum_{\ell=-\infty}^r$
instead of $\sum_{\ell=-\infty}^0$),
$$P_\lambda \{ S(B) \geq \lfloor \lambda c \rfloor - r \} = O(
\lambda^{-1/2} e^{-\lambda I}). \leqno({\rm A}.23)
$$
By adding up ({\rm A}.22), ({\rm A}.23) and applying Lemma 1, we can conclude that
$K_1 < \infty$. Then by (3.12), ({\rm A}.3) and ({\rm A}.4), $K \leq K_1 < \infty$.

Next, select $\epsilon$ small enough such that $\delta :=
[2 \pi \sigma_d(B) M''(\theta_c)]^{-1/2} \sum_{\ell=-\infty}^0 
e^{\theta_c \ell} - 
\epsilon > 0$. By Lemma A.2 with $r=0$ and 
$L > \sup_{\bx,\by \in D} \| \bx-
\by \|_\infty$, there exists $k$ large enough
such that
$$\sum_{\bv \in (k \lambda^{-1} {\bf Z})^d: 0 < \| \bv \|_\infty
\leq L} P_\lambda \{ S(B) \geq \lfloor \lambda c \rfloor,
S(\bv+B) \geq \lfloor \lambda c \rfloor \} \leq \epsilon
\lambda^{-1/2} e^{-\lambda I} \leqno({\rm A}.24)
$$
for all large $\lambda$. Then by stationarity, (3.10) and ({\rm A}.24),
noting that $x_\lambda \geq 0$,
$$\begin{aligned} P_\lambda(E_{{\bf 0},m,\lambda}) \geq &
P_\lambda \{ S(\bu+B) \geq \lfloor \lambda c \rfloor \mbox{ for
some } \bu \in (k \lambda^{-1} {\mZ})^d, {\bf 0} \prec \bu \prec
m \lambda^{-1} {\bf 1} \} \cr
\geq & \sum_{\bu \in (k \lambda^{-1} {\mZ})^d: {\bf 0} \prec
\bu \prec m \lambda^{-1} {\bf 1}} \Big( P_\lambda \{ S(\bu+B) \geq
\lfloor \lambda c \rfloor \} \cr
& \qquad - \sum_{\bw \in (k \lambda^{-1} {\mZ})^d:
0 < \| \bw-\bu \|_\infty \leq L} P \{ S(\bu+B) \geq \lfloor \lambda
c \rfloor, S(\bw+B) \geq \lfloor \lambda c \rfloor \} \Big) \cr
\geq & (m/k)^d \delta \lambda^{-1/2} e^{-\lambda I}
\mbox{ for all large } \lambda,
\end{aligned}
$$
and by letting $m \rightarrow \infty$ with $k$
fixed, it follows from (3.3), (3.11) and $x_\lambda$
bounded that $K > 0$. $\Box$

\bigskip
\centerline{\sc References}
\medskip

\ref {\sc Adler, R.} (2000). On excursion sets,
tube formulas and maxima of random field. {\it
Ann. Appl. Probab.} {\bf 10} 1-74.

\ref {\sc Albin, J.M.P.} (1990). On extremal
theory for stationary processes. {\it Ann.
Probab.} {\bf 18} 92-128.

\ref {\sc Aldous, D.} (1989). {\it Probability Approximations
via the Poisson clumping heuristic}. Springer-Verlag, New York. 

\ref {\sc Alm, S.E.} (1997). On the distribution
of scan statistics of a two-dimensional Poisson
process. {\it Adv. Appl. Probab.} {\bf 29}
1-18.

\ref {\sc Anderson, N.H.} and {\sc Titterington,
D.M.} (1997). Some methods for investigating
spatial clustering with epidemiological applications.
{\it JRSS `A'} {\bf 160} 87-105.

\ref {\sc Berman, S.M.} (1982). Sojourns and
extremes of stationary processes. {\it Ann.
Probab.} {\bf 10} 1-46.

\ref {\sc Bickel, P.} and {\sc Rosenblatt, M.} (1973).
Two-dimensional random field. In {\it Multivariate
Analysis III} (P.R. Krishnaiah, ed.) 3-15. Academic,
New York.

\ref {\sc Chan, H.P.} and {\sc Lai, T.L.} (2006).
Maxima of asymptotically Gaussian random fields
and moderate deviation approximations to 
boundary-crossing probabilities of sums of
random variables with multidimensional indices.
{\it Ann. Probab.} {\bf 34} 80-121.

\ref {\sc Chan, H.P.} and {\sc Zhang, N.R.} (2007).
Scan statistics with weighted observations. {\it Jour. Amer.
Statist. Assoc.}, {\bf 102} 595-602.

  \ref {\sc Cressie, N.} (1993). {\it Statistics
for Spatial Data}. Wiley, New York.

\ref {\sc Feller, W.} (1971). {\it An Introduction
to Probability Theory and its Applications}, Vol 2.
Wiley, New York.

\ref {\sc Glaz, J., Naus, J.} and {\sc Wallenstein, S.}
(2001). {\it Scan Statistics}. Springer, New York.

\ref {\sc Hogan, M.} and {\sc Siegmund, D.O.} (1986).
Large deviations for the maxima of some random fields. 
{\it Adv. Appl. Math.} {\bf 7} 2-22.

\ref {\sc Loader, C.} (1991). Large-deviation
approximation to the distribution of scan statistics.
{\it Adv. Appl. Probab.} {\bf 23} 751-771.

\ref {\sc Naus, J.I.} (1965). Clustering of random
points in two dimensions. {\it Biometrika} {\bf 52}
263-267.

\ref {\sc Pickands, J.} (1969). Upcrossing probabilities
for stationary Gaussian processes. {\it Trans. Amer.
Math. Soc.} {\bf 145} 51-75.

\ref {\sc Piterbarg, V.} (1996). {\it Asymptotic 
Methods in the Theory of Gaussian Processes and Fields}.
Transl. vol 148, Amer. Math. Soc., Providence.

\ref {\sc Qualls, C.} and {\sc Watanabe, H.}
(1973). Asymptotic properties of Gaussian random
fields. {\it Trans. Amer. Math. Soc.} {\bf 177}
155-171.

\ref {\sc Shafie, K.}, {\sc Sigal, B.}, {\sc Siegmund,
D} and {\sc Worsley, K.J.} (2003). Rotation
space random fields with an application to fMRI data.
{\it Ann. Statist.} {\bf 31} 1732-1771.

\ref {\sc Siegmund, D.O.} (1985). {\it Sequential
Analysis}. Springer, New York.

\ref {\sc Siegmund, D.O.} and {\sc Worsley, K.J.} (1995).
Testing for a signal with unknown location and
scale in a stationary Gaussian random field.
{\it Ann. Statist.} {\bf 23} 608-639.

\ref {\sc Spivak, M.} (1965). {\it Calculus on
Manifold: A Modern Approach to Classical Theorem
of Advanced Calculus}. W.A. Benjamin, New York.

\ref {\sc Taylor, J.} (2006). A Gaussian kinematic formula.
{\it Ann. Probab.} {\bf 34} 122-158.

\ref {\sc Taylor, J.}, {\sc Takemura, A.} and 
{\sc Adler, R.} (2005). Validity of the expected
Euler characteristic heuristic. {\it Ann.
Probab.} {\bf 33} 1362-1396.

\end{document}